# Generalized Contingency Analysis Based on Graph Theory and Line Outage Distribution Factor

Mohammad Rasoul Narimani, *Member, IEEE*, Hao Huang, *Student Member, IEEE,* Amarachi Umunnakwe, *Student Member,* Zeyu Mao, *Student Member*, Abhijeet Sahu, *Student Member,* Saman Zonouz, *Senior Member,* and Kate Davis, *Senior Member, IEEE*

***Abstract*-- Identifying the multiple critical components in power systems whose absence together has severe impact on system performance is a crucial problem for power systems known as (N-x) contingency analysis. However, the inherent combinatorial feature of the N-x contingency analysis problem incurs by the increase of x in the (N-x) term, making the problem intractable for even relatively small test systems. We present a new framework for identifying the N-x contingencies that captures both topology and physics of the network. Graph theory provides many ways to measure power grid graphs, i.e. buses as nodes and lines as edges, allowing researchers to characterize system structure and optimize algorithms. This paper proposes a scalable approach based on the group betweenness centrality (GBC) concept that measures the impact of multiple components in the electric power grid as well as line outage distribution factors (LODFs) that find the lines whose loss has the highest impact on the power flow in the network. The proposed approach is a quick and efficient solution for identifying the most critical lines in power networks. The proposed approach is validated using various test cases, and results show that the proposed approach is able to quickly identify multiple contingencies that result in violations.**

***Index Terms*-- Contingency analysis, line outage distribution factors, graph theory, betweenness centrality.**

## I. Introduction

THIS work is motivated by extreme events in electric power systems that are caused by multiple contingencies. Robust operation of the power grid requires anticipation of unplanned component outages that could trigger extreme events [1].Although electric power grids are designed to be resilient against any single contingency (N-1 contingency), loss of multiple components can lead to system instability, uncontrolled separation, cascading outages, voltage collapse, etc [1].Additionally, smart grid deployment exposes the electric power grid to purposeful and malicious attacks which raises concerns about the possibility of malicious N-x scenarios. Thus, it is important to make the system secure not only for any given N-1contingency but also for selected N-x contingencies. This results in a large number of contingencies that need to be analyzed [2]. Numerous efforts have been put on identifying the most critical components in electric power grids. Deter-mining and evaluating all possible combinations of component failures is a combinatorial problem which is not tractable even for medium size power systems. This is easily demonstrated by direct calculation of the number of combinations required for N-3 and N-4 analysis for a 10000-component system. The N-3 analysis involves approximately $\frac{10^{12}}{6}$ combinations and the N-4 analysis involves $\frac{10^{16}}{24}$ [3]. In this paper we propose a tractable generalize contingency selection approach based on the graph theory concept, i.e. group betweenness centrality (GBS), together with line outage distribution factors that identifies group of components whose loss would have severe impact on the power system. We hereafter refer to those components as critical components.

Critical components can be identified through network structural analysis. Graph-based methodologies provide promising approaches capable of finding these critical components, e.g. lines, buses in power systems [4]. Graph topology can also be used to detect anomalies in electric power grid [5]. Previous work has proposed a variety of metrics to identify the most critical components in electric grid. Developed in the context of social network analysis, betweenness centrality is a concept that captures the relative importance of an entity in a network [6]. These include betweenness centrality that reflects the edge or node importance in network structure [7].

Multiple studies utilize the betweenness centrality concept to identify the most critical component in power grid [4], [8]. The proposed approach in [4] applies the network centrality measures to fulfill the N-1 contingency analysis. Similarly, the approach in [8] uses different centrality measure including node and edge betweenness centrality to identify important nodes and edges in power systems. Electrical centrality metric which is calculated based on the impedance matrix $Z_{bus}$ and utilizes the centrality metric, to explain why in electric power grids a few numbers of highly connected bus failures can cause cascading effect, is investigated in different studies [9]–[11]. A graph edge betweenness centrality measure is proposed in [12] to perform contingency analysis of large-scale power grid. The betweenness centrality concept is extended in [13] to account for N-x contingency selection for $x \geq 2$.

Although the structural analysis of electric power grid returns invaluable information regarding the vulnerability analysis, pure structural analysis cannot capture the physical features of the electric power grid. Thus, it is essential to have a more comprehensive model for approximating failure behavior of electric power grid [14]. To address these issues, different properties of electric power grids along with the betweenness

The work described in this paper was supported by funds from the US Department of Energy under award DE-OE0000895 and the National Science Foundation under Grant 1916142.

centrality concept have been taken into consideration to address the physics of the networks. Maximal load demand and the capacity of generators are considered along with the betweenness centrality in [15] to analyze the vulnerability of electric power grid. Electrical distance, an essential feature of power network, is taken into account in the betweenness centrality metric to capture the grid properties [16]. The betweenness centrality is applied to the weighted graph of electric power grid, weighted by the corresponding admittance matrix, for vulnerability analysis of power network [17].

Although the proposed approaches in [15]-[17] consider the electrical properties of electric power grids, none of them captures the impacts that loss of a component might have on the system. In this paper, the line outage distribution factors(LODF) [18], a sensitivity metric of how a change in a line's status affects the flows on other lines in the system, enables the proposed approach to capture the physics of the power network and thus improves the accuracy of the results. The LODF metric has already been used in identifying multiple contingencies in power systems [18], but the approach proposed in [18] is limited to N-2 contingency analysis. Additionally, in this paper the betweenness centrality factor is extended to the group betweenness centrality which facilitate searching for multiple components in the network. The proposed approach in this paper is generalized contingency analysis that without any limitation performs different multiple contingency analysis.

The proposed method in this paper is demonstrated for the line outages, but without loss of generality it can also be used in similar way for node outages. The contribution of this paper is threefold. First, we present a novel scheme to apply N-x group betweenness centrality to contingency analysis. Second, we leverage LODF metrics to facilitate the proposed algorithm with augmenting line outage impact on power flow in the network. Third, we validate that the proposed method is computationally tractable for multiple (N-x) contingency analysis in very large power systems, i.e. with couple of thousands of lines.

This paper is organized as follows. Section II reviews the betweenness centrality measure and its extension, i.e. group betweenness centrality. Section III describes the proposed generalized contingency analysis (GCA) underlying our contribution in finding the most critical lines in the network. Section IV empirically evaluates the proposed approach. Section V concludes the paper.

## II. CENTRALITY DEFINITION AND EXTENSION

Centrality metrics are used in network science to rank the relative importance of vertices and edges in a graph. In network analysis, there are several metrics for the centrality of a vertex or an edge [8]. In this paper we use group betweenness centrality of the edges to find the most critical lines in power systems. We first review the definition of betweenness centrality and then extend it to the group betweenness centrality.

### A. Betweenness Centrality

Electric power grids can be modeled as graphs G:=(V,E), where V(vertices) and E(edges) are the sets of buses and lines. Betweenness centrality measures the extent to which a vertex lies on paths between other vertices. Vertices with high betweenness may have considerable influence within a network since more paths that connect different vertices passthrough them. They are also the ones whose removal from the network will most likely disrupt information or power flows between other vertices because they lie on the largest number of paths connecting different vertices.

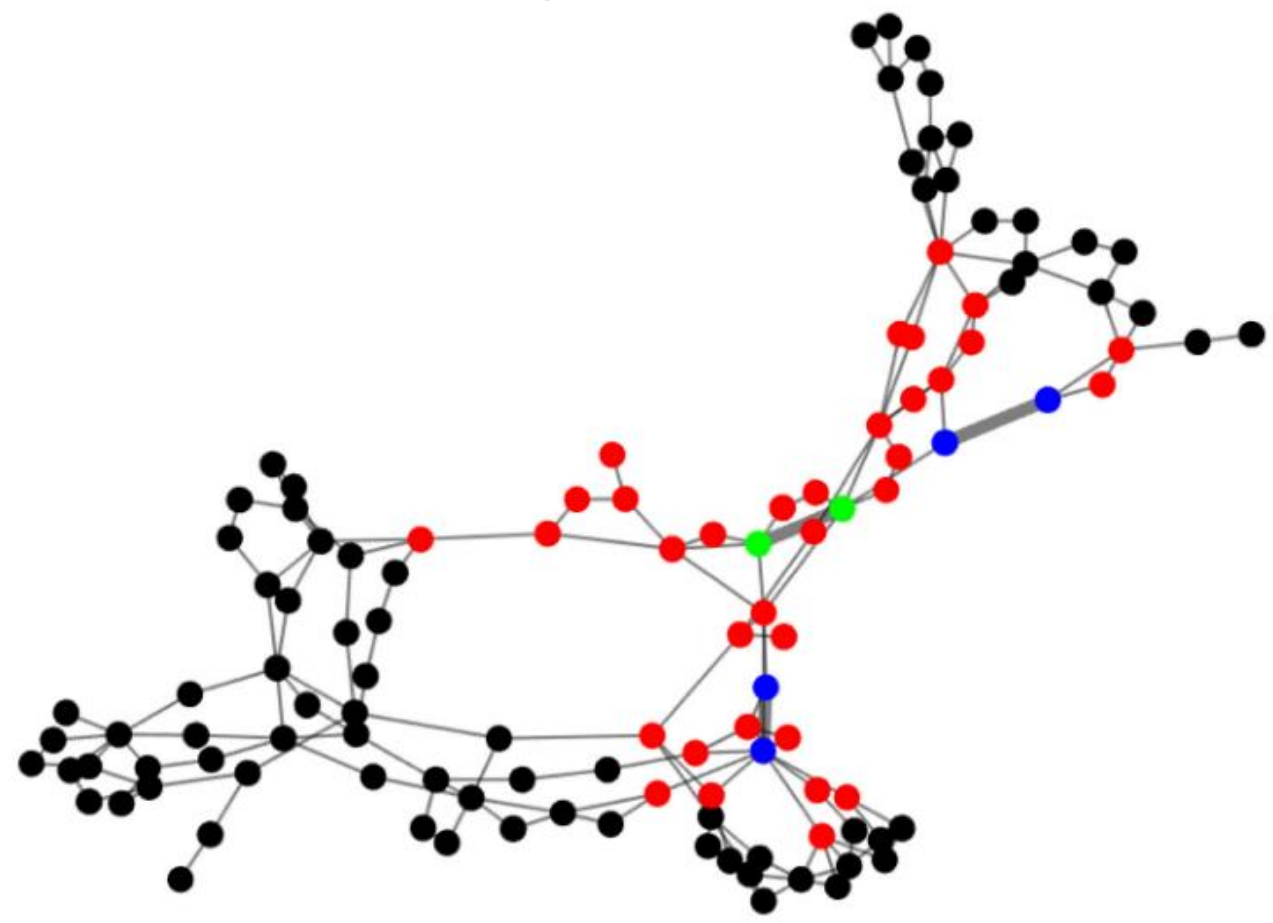

Figure 1. The equivalent graph of IEEE 118-bus test cases. The green nodes show the both ends of the line whose outage has the highest LODF factor. All the neighbors with three hop distance (i.e. search level equals to three) from the green nodes are shown by red colors. The lines with the blue end node s are the second and the third lines with the highest LODF values. These lines are used to implement the concept of the group betweenness centrality for the N-3 contingency analysis. The set of red nodes will be investigated to determine the critical lines.

The betweenness centrality concept can be extended from vertices to measure the influence of edges in a graph. Betweenness centrality is defined as the ratio of the number of shortest paths that passthrough an edge to the total number of shortest paths between all possible pairs of vertices. Figure 1 visualizes the concepts of the betweenness centrality and its expansion, i.e. group betweenness centrality, to makes them more understandable. Mathematically, betweenness centrality of an edge can be expressed as:

$$BC(e) = \sum_{\{s,t\in V\}} \frac{\sigma(s,t|e)}{\sigma(s,t)} \qquad (1)$$

where $\sigma(s,t|e)$ is the number of shortest paths in the graph between s and t that contain edge e, and $\sigma(s,t)$ represents the number of shortest paths in the graph between s and t.

The betweenness centrality concept measures the importance of a single node (edge) and thus can be applied to solve the N-1contingency analysis. The underlying assumption for the proposed N-x contingency analysis is that cascading failures might occur from the simultaneous failure of a x critical lines.

Thus, the betweenness centrality concept needs to be adapted to measure the importance of multiple nodes (edges). In this connection, we proposed a new group betweenness centrality that is explained in sequel.

### B. Group Betweenness Centrality

The goal of the group betweenness centrality is to identify a set of the most important components that their lost has a severe impact on the network. However, betweenness centrality metric is for the individual component (a vertex or an edge).

Thus, to perform the N-x contingency analysis we extend the betweenness centrality metric to consider a group of entities. Mathematically, group betweenness centrality of an edge can be expressed as:

$$GBC\left(E_G^i\right) = \sum_{\{s,t \in V \setminus E_G^i\}} \frac{\sigma\left(s,t\middle|E_G^i\right)}{\sigma(s,t)} \qquad (2)$$

where $E_G^i$ is a subset of edges of interest, $\sigma(s,t)$ is the number of shortest paths between s and t, and $\sigma\left(s,t\middle|E_G^i\right)$ is the number of shortest paths between s and t that contain any element in $E_G^i$. The notion of group betweenness centrality first introduced in [19], [20] to identify groups of individuals who have collective influence in a social network. For implementing the group betweenness centrality in (2) we first need to determine groups of lines, with the same number of lines (i.e. x), that need to be evaluated. However, the number of groups that need to be evaluated increases exponentially with respect to the number of elements in the group, $\binom{N}{x}$, which makes the group betweenness centrality in (2) computationally intractable for the large test cases. To cope with this issue, we leverage the line outage distribution factor, explained in section III, to select limited numbers of groups out of the astronomical group of components that needs to be evaluated. Leveraging LODF to select group of lines that their loss has more impact on transmission power systems drastically reduces the number of groups that need to be evaluated which is elaborated in Section IV.

## III. SOLUTIONMETHODOLOGY

This section leverages the group betweenness centrality concept presented in Section II-B and the LODF metric to identify the most critical lines in the electric power grids. The pure graph theory information, i.e. group betweenness centrality, cannot fully address the characteristic of the electric power grid. Thus, it is essential to incorporate the LODF metric to take the electric power grid's features into account. The next section reviews the LODF metric first and then presents the proposed algorithm that leverages LODF and Group Betweenness Centrality to select critical N-x contingencies.

### A. Line outage distribution factors

Line Outage Distribution Factors (LODFs) are a sensitivity measure of how a change in a line's status affects the flows on other lines in the system. LODFs are used extensively when modeling the linear impact of contingencies in PowerWorld Simulator [21]. LODF metrics in electric power grids provide approximate but quick solutions for the change in the line flow. The quick computation of the LODF metric makes it an attractive measure for solving different problems in power systems. The LODF metric is used to screen multiple element contingencies in [18]. The LODF is also used for detecting island formation in power networks [22]. The LODF used to solve the security constrained unit commitment problem [23]and for evaluating distribution network expansion options [24].

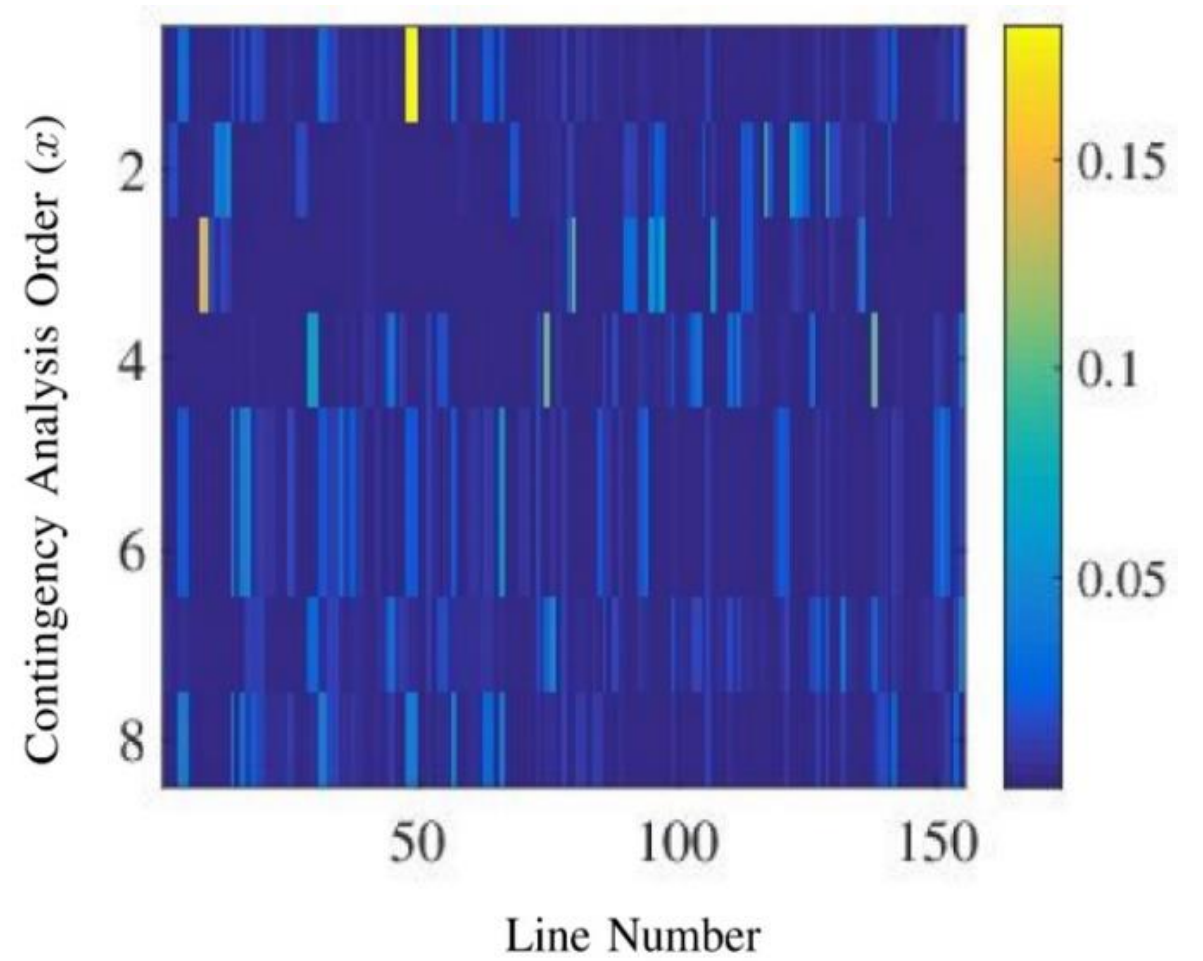

Figure 2. Normalized LODF for different order of contingency analysis (i.e. different x in the N-x term) for 200-bus test cases. The colors demonstrate the normalized LODF values.

LODFs vary with change in the topology, when an outage occurs [25]. Multiple efforts have studied calculating LODFs dynamically, i.e. after occurring outages, as well as extending the LODF definition for multiple contingency analysis [25]-[27]. Although the LODF metric changes after outages, these changes usually happen gradually as the number of outages increase. Therefore, it would be a reasonable assumption for the LODF metric to remain fixed for a few line outages. To validate this assumption, we investigate the LODF changes for different N-x contingency analysis for 200-bustest systems which are illustrate in Figure 2. The x- and y-axis show the line number and the contingency analysis order (i.e. x in the *N-x* term), respectively. The colors demonstrate the normalized LODF values (i.e. LODF values divided by maximum value of LODFs after each outage). Lines with very small LODFs are mitigated for illustration purposes. Figure 2shows that the LODF would not change drastically after few line outages which validates our assumption.

### B. Proposed Methodology

We leverage the group betweenness centrality concept and the LODF metric to identify the most critical lines in contingency analysis. The group betweenness centrality has previously been proposed to find critical lines in a power system [13]. However, the procedure in [13] uses pure topological information to find critical lines but neglects the electrical characteristics. Moreover, the procedure utilized to identify the groups is complicated. These two shortcomings are addressed in this paper via incorporating the LODF metric into our GBC methodology. To this end, we first find the LODF metric for all lines and then select the limited number of lines based on the measure in (3). To further capture the physics of power networks and prevent selecting lines that carry a small portion of power flow, the power flows in lines are incorporated in (3).The mean of the remaining lines LODFs seems to be an appropriate measure for determining the importance of the lost line in the system. However, the value of LODF can be either
positive or negative, so the mean of the corresponding LODFs may not capture the importance due to the offset of positive and negative LODF values for different lines. Therefore, in this paper, we use the absolute value of LODFs to capture the line importance in power systems. Furthermore, losing a line might change the LODF values of a small set of the lines drastically while change the others slightly. We are looking fora line that its lost increases the LODF of all other lines not only a small set of them. To this end, the standard deviation of the absolute LODFs values is incorporated in the proposed metric for identifying the high impacted lines in the system. Mathematically, the proposed metric can be expressed as:

$$NLODF(i) = \frac{mean(abs(LODFs))}{std(abs(LODFs))} \qquad (3)$$

$$M(i) = PF(i) \times \min\{NLODF(i), 1\} \qquad (4)$$

Where NLODF and M are normalized LODF metric and the proposed measure for selecting critical lines, respectively. $PF(i)$ is the power flow in line $i$ during the normal operation; $mean(.)$, $std(.)$, and $abs(.)$ indicate the mean, the standard deviation, and the absolute functions, respectively. In this paper, we utilize PowerWorld Simulator [21] to calculate LODFs. However, if losing a line creates an islanding situation, the PowerWorld Simulator does not calculate the exact LODF but assign a large value to the remaining lines' LODF which makes the $NLODF(i)$ value a very large number. Without islanding situation, the $NLODF(i)$ values varies between [0,1]. For handling the islanding situation, we enforce the $NLODF(i)$ value corresponds to the islanding situation to be 1 as it is formulated in (3b). Equation (3b) represents the metric that is used to identify critical lines.

A metric that leverages the physics of the power networks as well as their topological information can effectively identify the critical lines in power systems. To this end, we first evaluate equation (3b) for all the lines in the network and select the first 10% of the lines with the highest value. Then we search for the critical lines based on the search level i.e. a prespecified parameter that determine how deep the proposed algorithm should search around the selected lines.

Figure 1 visualizes the searching mechanism of the proposed approach for the IEEE 118-bus test case. We first select the lines with highest measure value of (3b). One of these lines is illustrated with green buses at its both ends. The next step is finding all the nodes in graph G that are within distanced, i.e. search level, from these green nodes. For brevity this set of buses are named "neighboring buses" in this paper. For the N-1contingency analysis, equation (1) is applied to the neighboring buses where *s*, *t* are different combination of buses in this set and e is the line that both ends are shown by the green buses. For the N-x contingency analysis we should first determine the $E_G^i$, the subset of edges of interest. In this connection, other lines with high value of (3b) in which their end buses are within the neighboring buses are selected for the N-x contingency analysis. It is notable that there are two parameters that we can control to find the critical lines in the proposed approach. The first parameter is x which determine the contingency analysis level. The second parameter is the "search level" that determines the number of neighboring buses that should be evaluated. The higher search level increases the search distance in the graph for the contingency analysis, i.e. more red buses in Figure 1, but it incurs more computational burden. The execution time for different search level is elaborated for different test cases in Section IV.

**Algorithm 1** Discover Critical Branches in each sub graph

```
1. Initiate esa object and load case.
2. Extract Br= branch info, pf=branch power flow
3. Create Multigraph G(Br) = (V,E).
4. LODF = CalculateLODF(case).
for b in Br do
    LODFmetric_b = mean(LODF_b)/sd(LODF_b)
    if LODFmetric_b ⩾ 10 then
        LODFmetric_b = abs(pf_b)
    else
        LODFmetric_b = LODFmetric_b * abs(pf_b)
    end if
end for
5. SortedBranch = Sort b by LODFmetric_b
6. InvgtdBranch=K% branches from SortedBranch
7. Initialize Invgtd_Br_Nbr
for from, to, c_id, w in InvgtdBranch do
    from_sg = subgraph(G, from, search_level)
    to_sg = subgraph(G, to, search_level)
    Inv_Br_Nbr = from_sg.nodes ∪ to_sg.nodes
    Invgtd_Br_Nbr.append(Inv_Br_Nbr)
end for
7. crit_lines_in_sg =Importance_Subgraph(Invgtd_Br_Nbr)
8. return crit_lines_in_sg
```

In a cascading power transmission outage, component outages might propagate nonlocally; after one component outages, the next failure may be far, both topologically and geographically [28]. It is notable that the higher "search level" parameter enables the proposed approach to find critical lines that are far away from each other. However, higher "search level" incurs more computational complexity to the proposed approach especially when it applies on larger test cases. Another feature that enables the proposed approach to globally search for the critical lines is that it selects the first Y% of the lines (e.g. 10%) with higher measure value of (3b) and then starts searching around those lines. Those selected lines are usually distant, both topologically and geographically, that facilitates the proposed approach to find critical lines globally.

The steps of the proposed approach for identify critical lines are detailed in Algorithm 1 and 2.

**Algorithm 2** Importance_Subgraph ($Investigated_Br_Nbr$)

```
Initialize imp_line_crit_sg_list
for ibn in Investigated_Br_Nbr do
    inv_sg=subgraph(ibn)
    Initialize inv_sg_res dict
    for edge in inv_sg.edges do
        inv_sg_res(edge) = LODF(edge)
    end for
    sorted_inv_sg_res = sort(inv_sg_res)
    Initialize crit_lines array.
    Initialize x in N − x contingency
    for key, val in sorted_inv_sg_res[1 : x] do
        from = val.from
        to = val.to
        crit_lines.append(from, to)
    end for
    imp_line_crit_sg_cnt = 0
    for line in crit_lines do
        imp_line_cnt = 0
        for from_n in inv_sg.nodes do
            for to_n in inv_sg.nodes do
                sp = srtst_path(inv_sg, from_n, to_n)
                if line in sp then
                    imp_line_cnt+ = 1
                end if
            end for
        end for
        imp_line_crit_sg_cnt+ = imp_line_cnt/2
    end for
    imp_line_crit_sg_list.append(imp_line_crit_sg_cnt)
end for
return imp_line_crit_sg_list
```

## IV. RESULTS ANDANALYSIS

This section demonstrates the effectiveness of the proposed approach using selected synthetic test cases from the benchmark library for electric power grids in Texas A&M University [29]. These test cases were selected since they mimic realistic electric power grids. Synthetic electric grid cases are a representation of fictitious power grids with a detailed modeling of the power system elements [30]. Our implementations use PowerWorld [21] for contingency analysis, Python as a programming language, and ESA [31], a Python package that provides an easy to use and light-weight wrapper for interfacing with PowerWorld's Simulator Automation Server (SimAuto), as an interference to communicate with PowerWorld. The results are computed using a laptop with an i7 1.80 GHz processor and 16 GB of RAM. Three different test cases including 200-, 500-, and 2000-bustest cases are investigated to authenticate the ability of the proposed approach in finding critical lines in electric power grids. The result for each test case is discussed in detail in sequel.

- **200-Bus Test System**

The 200-bus test system is a synthetic test case with 245branches and 49 generators that builds from public information and a statistical analysis of real power systems [29]. This test system is selected purposely since its relatively small size enables the exhaustive search in identifying the critical lines which can be used for authenticating the results obtained by the proposed approach. In this connection, a brute-force search is performed for finding the N-1contingencies in the 200-bus test system to evaluate the performance of the proposed algorithm. It is notable that the 200-bus test system is N-1resilient and as it was expected the brute-force search finds only one reserved limit violation, where there is not enough active power reserves in the make-up power specification to cover the active power changes by the contingency. The application of the proposed approach finds the exact same violation for the N-1contingency analysis. This verifies the ability of the proposed approach in identifying contingencies in power systems. To further evaluate the ability of the proposed approach in finding the critical lines, another brute-force search is fulfilled for identifying N-2 contingencies in 200-bus test system and the obtained results match those obtained by the proposed algorithm. The brute-force search and the proposed approach found the N-2 contingencies for 200-bus test cases in 230 and 38 seconds, respectively. Comparing the computational times of the brute force search and the proposed approach for finding N-2 contingencies reveals the ability of the proposed approach in finding critical lines in large test cases which cannot be done by the exhaustive search. The proposed approach is applied to identify critical lines for different contingency analysis in the 200-bus test system. The results are summarized in Table I. The first column lists the x in the N-x contingency term. The second and the third columns represent critical lines and contingency types, respectively and the fourth column represents number of contingencies.

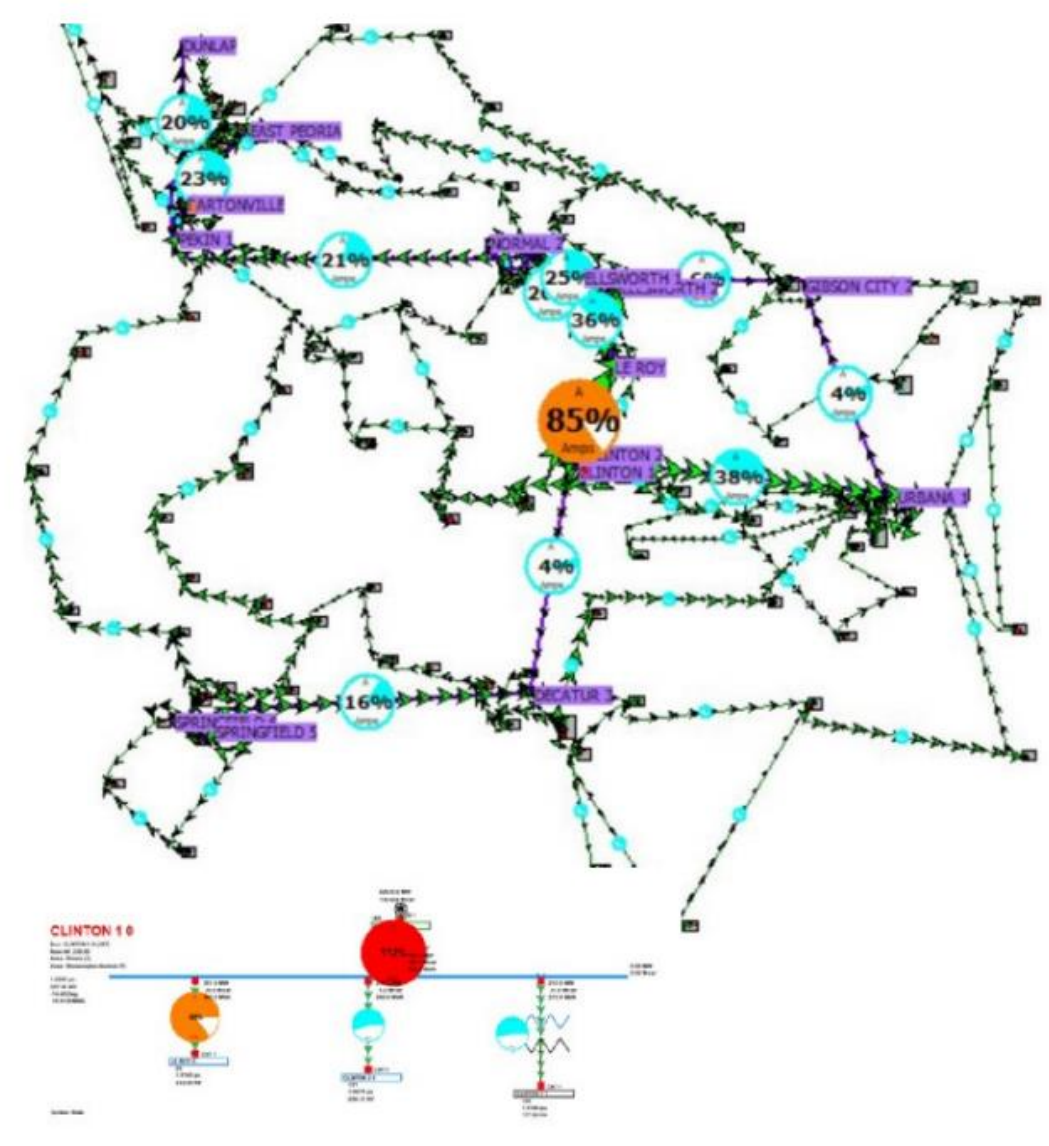

Figure 3. One-line diagram of the 200-bus test case and the corresponding overflow contingency caused by the removal of [136, 133], [135, 133] lines. The bottom portion of Figure 3 shows the zoom-in view of the area that violation has occurred.

With the contingency analysis tool in PowerWorld, we capture four types of "limit violations" in this paper, including reserve limit, overflow, undervoltage and unsolved. The unsolved cases represent a situation where there is no solution for the power flow equations, or they cannot converge. The overflow and undervoltage violations can be counted by the number of components that fall into the category. The critical lines in the second column are those that their lost has severe impact on the network performance. These data are crucial for power system operation and planning. The One-line diagram of the 200-bus test cases and the corresponding violations caused by the outage of lines [189, 187], [187, 121] are illustrated in Figure 3. The bottom portion of this figure shows the zoom-in view of the area that violation has occurred.

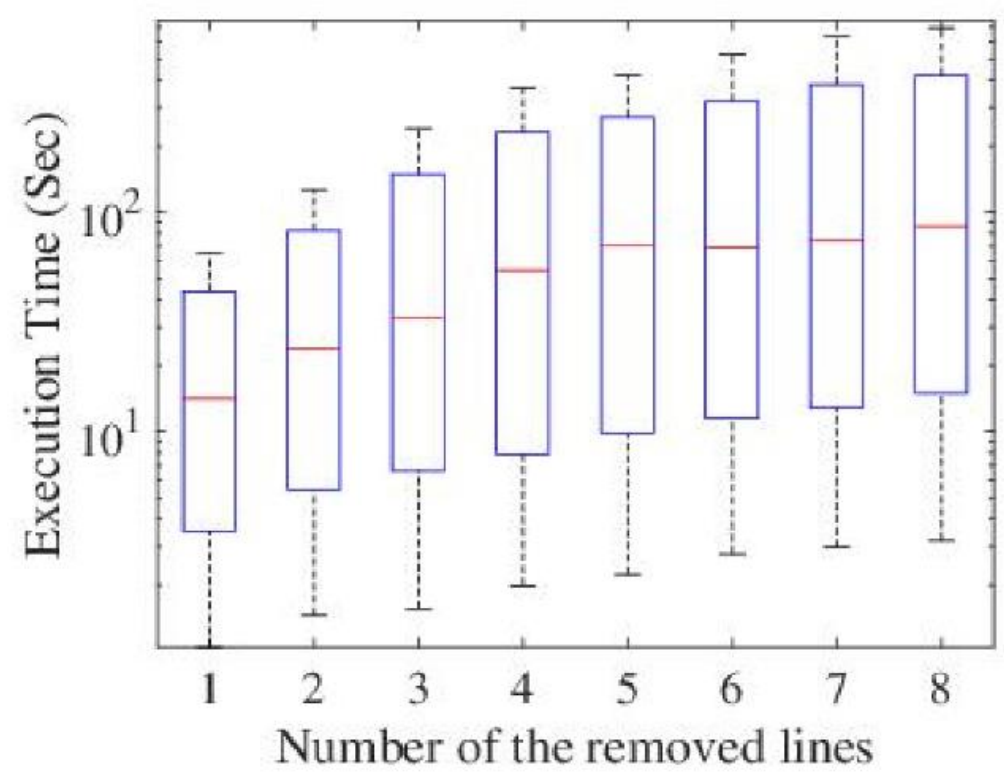

Figure 4. Execution time comparison for various x in the N-x contingency analysis of 200-bus test system.

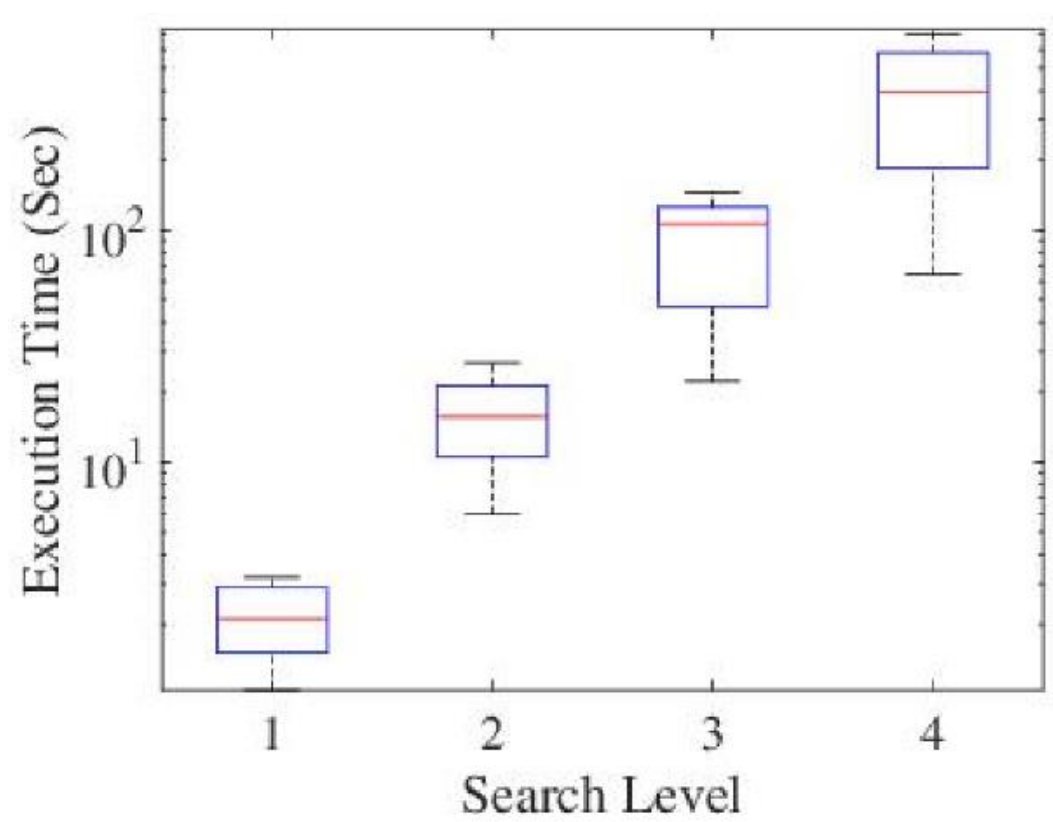

Figure 5. Execution time comparison for various “search level” in the N-x contingency analysis of 200-bus test system.

The box-and-whisker plot shown in Figure 4 compares the impact of x on the execution times of N-x contingency analysis for 200-bus test case. Every box contains the execution time of various contingency analysis with different search levels (i.e. 1-8) and a specific x in the N-x term. The lower and upper ends of the boxes in Figure 5 reflect the first and third quartiles, and the lines inside the boxes denote the median. There is not a significant difference between the median lines of the boxes which reflects the fact that the proposed approach can solve higher order of N-x contingency analysis within a reasonable time. This result was expected because for a same search level a same set of neighboring buses are utilized in identifying different N-x contingency analysis. In other words, the execution time of different N-x contingency analysis for specific search level changes slightly for different x values. The plot in Figure 5 shows the impact of different search levels on the execution times of contingency analysis in 200-bus test system. The execution time linearly increases with the search level increment. The increase in the search level enables the proposed approach to search in a broader area but it incurs computational burden.

Table I
RESULTS FROM APPLYING THE PROPOSED APPROACH TO 200-BUS TEST SYSTEM.

| x | Critical Lines | Violation Type | Number of Violations |
|---|---|---|---|
| 1 | [189, 187] | Reserve limit | NA |
| 2 | [189, 187],[187, 121] | Reserve limit | NA |
| 2 | [136, 133],[135, 133] | Overflow | 1 |
| 3 | [189, 187],[187, 121],[188, 187] | Reserve Limit | NA |
| 3 | [189, 187],[187, 121],[154, 149] | Reserve Limit | NA |
| 3 | [136, 133],[135, 133],[125, 123] | Overflow<br>Undervoltage | 2<br>18 |
| 4 | [189, 187],[187, 121],[188, 187],[179, 178] | Undervoltage | 25 |
| 4 | [189, 187],[187, 121],[154, 149],[152, 149] | Overflow | 2 |
| 4 | [136, 133],[135, 133],[125, 123],[126, 123] | Overflow<br>Undervoltage | 2<br>84 |
| 5 | [189, 187],[187, 121],[188, 187],[179, 178],[45, 187] | Undervoltage | 29 |
| 5 | [189, 187],[187, 121],[154, 149],[152, 149],[153, 149] | Unsolved | NA |
| 5 | [136,133],[135, 133],[125, 123],[126, 123],[127, 123] | Unsolved | NA |
| 5 | [125, 123],[126, 123],[127, 123],[124, 123],[134, 133] | Overflow<br>Undervoltage | 1<br>17 |
| 6 | [189, 187],[187, 121],[188, 187],[179, 178],[45, 187], [188, 89] | Undervoltage | 26 |
| 6 | [136,133],[135, 133],[125, 123],[126, 123],[127, 123],[188, 89] | Unsolved | NA |
| 6 | [189, 187],[187, 121],[154, 149],[152, 149],[153, 149],[155, 149] | Unsolved | NA |
| 7 | [136, 133],[135, 133],[125, 123],[126, 123],[127, 123],[124, 123],[134, 133] | Overflow<br>Undervoltage | 2<br>183 |
| 7 | [189, 187],[187, 121],[188, 187],[179, 178],[45, 187],[188, 89],[121, 178] | Undervoltage | 36 |
| 7 | [189, 187],[187, 121],[154, 149],[152, 149],[153, 149],[155, 149],[16, 15] | Unsolved | NA |
| 8 | [136, 133],[135, 133],[125, 123],[126, 123],[127, 123],[124, 123],[134, 133],[128, 133] | Unsolved | NA |
| 8 | [136, 133],[135, 133],[125, 123],[126, 123],[127, 123],[124, 123],[134, 133],[30, 29] | Overflow<br>Undervoltage | 2<br>151 |
| 8 | [189, 187],[187, 121],[154, 149],[152, 149],[153, 149],[155, 149],[16, 15],[188,187] | Unsolved | NA |
| 8 | [189, 187],[187, 121],[188,187],[45, 187],[188,89],[121,178],[14, 15],[188,187] | Unsolved | NA |

- **500-Bus Test System**

The 500-bus test system, with 597 branches and 90 generators, is a relatively large synthetic test case that is used for evaluating the ability of the proposed approach in identifying critical lines. All the synthetic test cases used in this paper are N-1 resilient, i.e. no contingency occurs with the outage of a single line. Furthermore, these test cases are resilient enough that it is hard to identify a contingency by randomly removing multiple lines. However, the proposed algorithm can identify overflow violations which occur by the outage of three branches, i.e. [142, 141], [424, 423], [87, 141]. Finding contingencies caused by the outage of three lines reveals the ability of the proposed approach in identifying critical lines in relatively large resilient test cases such as the 500-bus test system. The results of different contingency analysis for finding the most critical lines in the 500-bus test system are summarized in Table II. It is notable that multiple critical lines are identified for each N-x contingency analysis. The results in Table II summarize both contingency type and number of contingencies for the outage of critical lines in the 500-bus test system.

Table II
RESULTS FROM APPLYING THE PROPOSED APPROACH TO 500-BUS TEST SYSTEM.

| x | Critical Lines | Violation Type | Number of Violations |
|---|---|---|---|
| 3 | [142, 141],[424, 423],[87, 141] | Overflow | 3 |
| 4 | [162, 220],[23, 386],[87, 141], [247, 246] | Unsolved | NA |
| 5 | [162, 220],[23, 386],[87, 141], [247, 246],[437, 428] | Unsolved | NA |
| 6 | [197, 195],[198, 195],[162, 220], [456, 453],[23, 386],[247, 246] | Overflow | 5 |
| 6 | [162, 220],[23, 386],[87, 141], [247, 246],[437, 428],[438, 428] | Unsolved | NA |
| 7 | [197, 195],[198, 195],[162, 220], [456, 453],[23, 386],[247, 246], [164, 163] | Overflow | 5 |
| 7 | [162, 220],[23, 386],[87, 141], [247, 246],[437, 428],[438, 428], [439, 428] | Unsolved | NA |
| 8 | [197, 195],[198, 195],[162, 220], [456, 453],[23, 386], [247, 246], [164, 163],[455, 453] | Overflow | 5 |

The execution time of the contingency analysis approaches plays an important role in power system operation and control. Many approaches for identifying critical lines are not tractable to be applied on large test cases. However, the proposed approach performs the N-3 contingency analysis in the relatively large 500-bus test case in just 5 minutes.

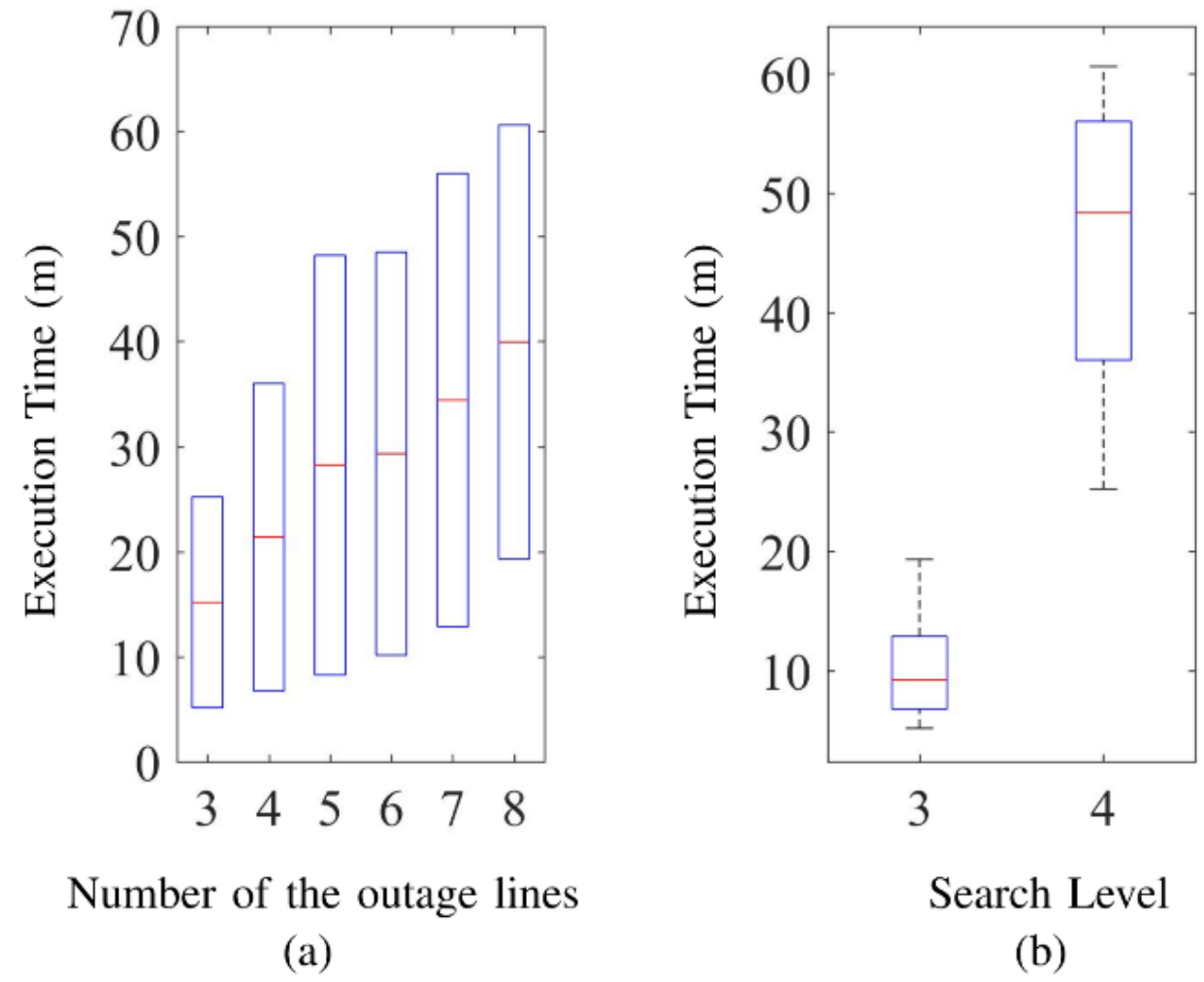

Figure 6. (a) Execution time comparisons of various contingency analysis using GCA. (b) Execution time comparisons of different search levels in GCA for 500-bus test system.

The One-line diagram of the 500-bus test cases and the corresponding violations caused by the outage of the lines [142, 141], [424,423], [87, 141] are illustrated in Figure 7. The bottom portion of this Figure 7 depicts the zoom-in view of the area that violation has happened ( i.e. the 112% overflow in the line).e plot in Figure 6 (a) compares the impact of x on the execution time of N-x contingency analysis for 500-bustest case. As it was expected the execution time of N-x contingency analysis for different x values does not change significantly for a same "search level". Similarly, the plot in Figure 6 (b) shows the impact of different search levels on the execution time of N-x contingency analysis for 500-bus test system. Although the increment in search level increases the execution time of N-x contingency analysis but the proposed algorithm remains computationally tractable for higher search levels.

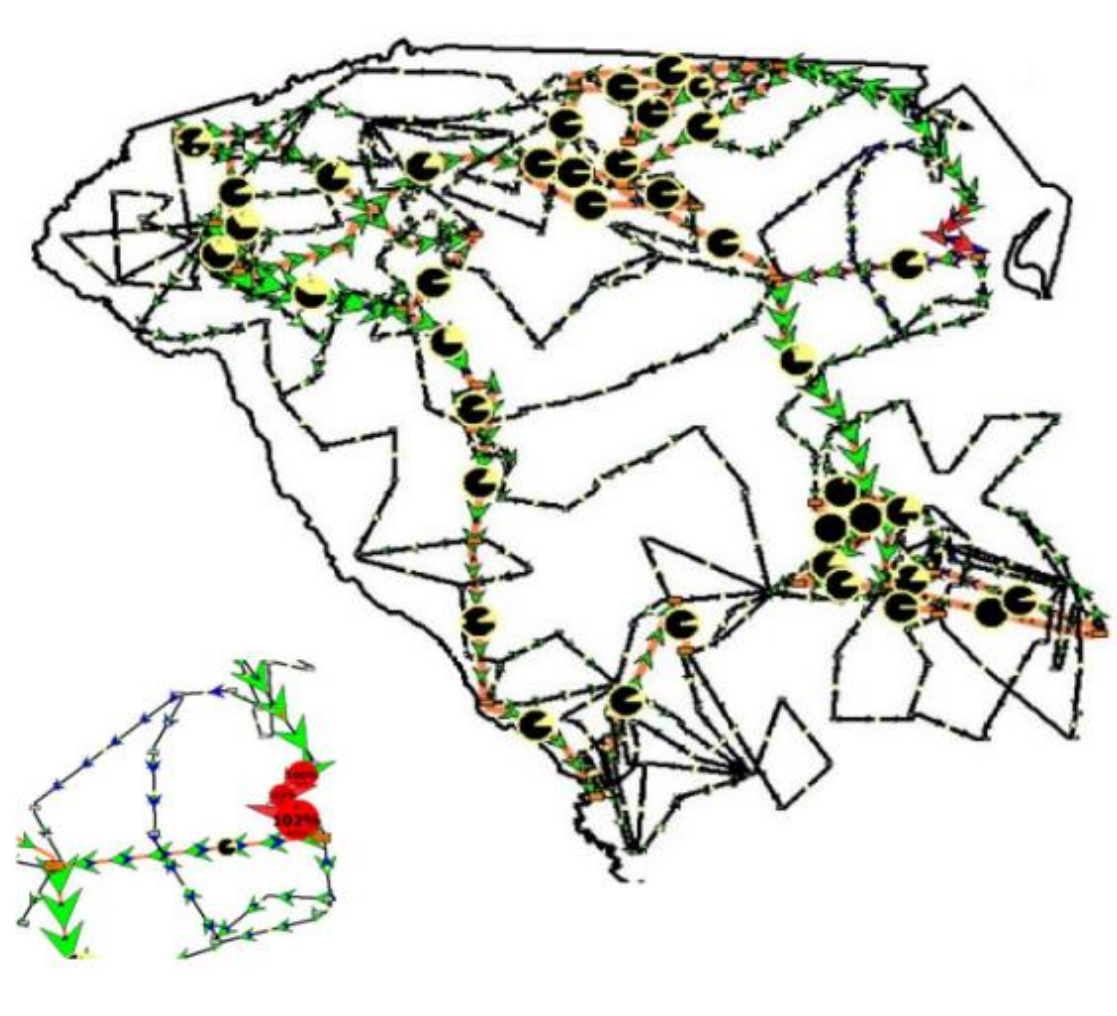

Figure 7. One-line diagram of the 500-bus test cases and the corresponding violations caused by the removal of [142, 141],[87, 141],[424, 423] lines. The bottom portion of Figure 7 shows the zoom-in view of the area that violation has occurred.

- **2000-Bus Test System**

The 2000-bus test system, with 3206 branches and 544 generators, is a synthetic test case that is utilized for evaluating the ability of the proposed approach in identifying critical lines in large test systems. The size of this synthetic test system makes it a challenging test case for contingency analysis, because the number of combinatorial scenarios that need to be evaluated increases drastically by the size of the network. Thus, for this test case the exhaustive search approaches are not tractable even for the N-3 contingency analysis. Like other synthetic test cases investigated in this paper, the 2000-bus test case is N-1 resilient. Furthermore, this test system is resilient enough that finding a contingency even after multiple lines outage, e.g. randomly removing 7 lines (N-7contingency), is not easy. Conversely, the proposed approach can easily find a N-4 contingency caused by the outages of four lines including [5262, 5260], [5263, 5260], [5317, 5260], [5358, 5179]. The impact of x and search level on the execution time of N-x contingency analysis for 2000-bus test case are shown in Figure 8 (a) and Figure 8 (b), respectively. For the same search level, execution time of N-x contingency analysis in Figure 8 (a) change linearly for different x values. This enables the proposed approach to evaluate higher orders of N-x contingency analysis in large test systems. The lower search level mirrors the lower number of lines that need to be investigated. Thus, for different search levels in Figure 8 (b), the lower search level value yields the faster N-x contingency analysis.

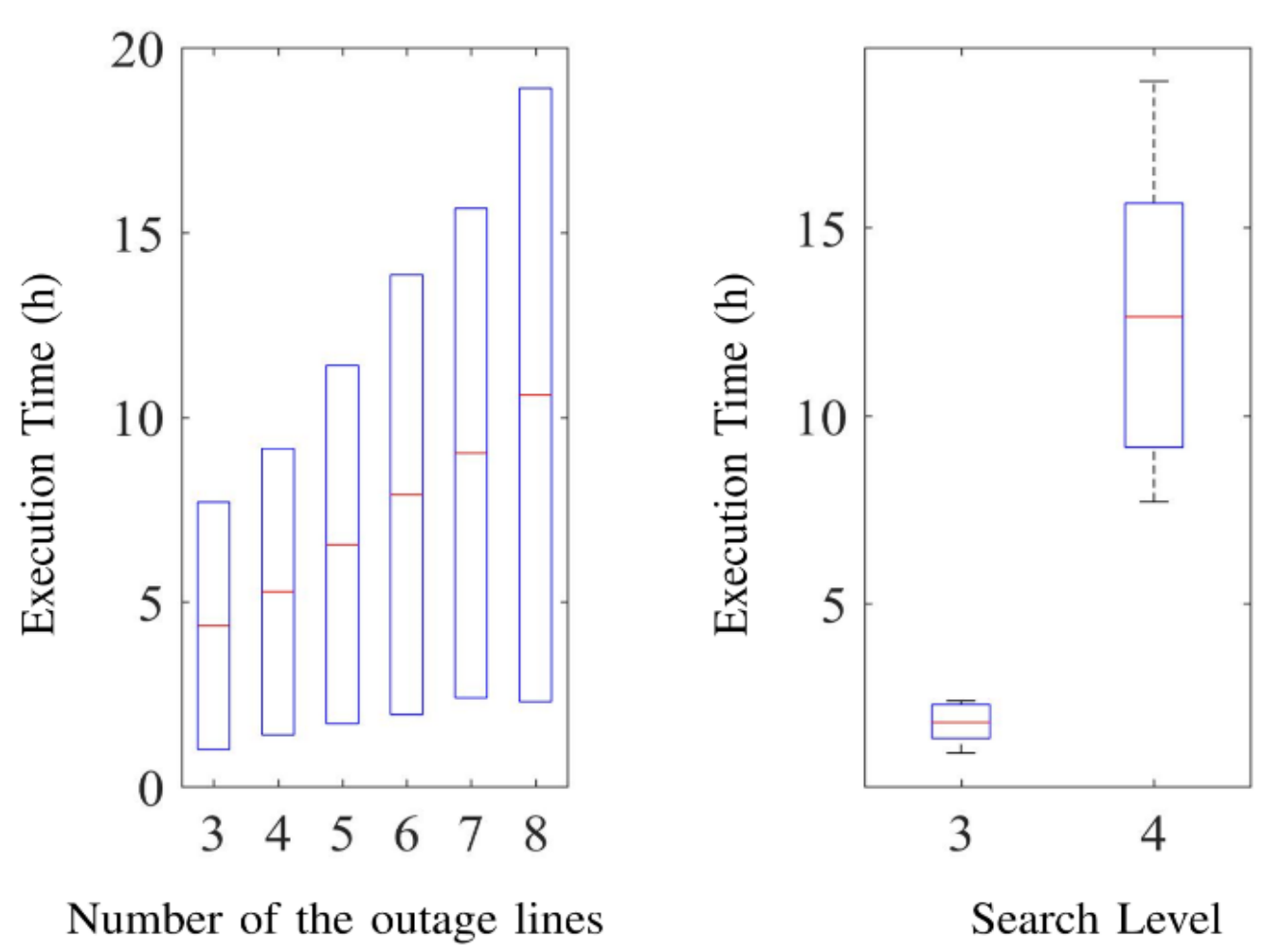

Figure 8. (a) Execution time comparisons of various contingency analysis using GCA. (b) Execution time comparisons of different search levels in GCA for 2000-bus test system.

Table III
RESULTS FROM APPLYING THE PROPOSED APPROACH TO 2000-BUS TEST SYSTEM.

| x | Critical Lines | Violation Type | Number of Violations |
|---|---|---|---|
| 3 | [7099, 7095],[7098, 7095],[7058, 7095] | Reserve Limit | NA |
| 4 | [5262, 5260],[5263, 5260],[5317, 5260], [5358, 5179] | Overflow | 4 |
| 4 | [7099, 7095],[7098, 7095],[7058, 7095], [6161, 7018] | Reserve Limit | NA |
| 4 | [7099, 7095],[7098, 7095],[7058, 7095], [7311, 7304] | Reserve Limit | NA |
| 5 | [5262, 5260], [5263, 5260], [5317, 5260], [5358, 5179], [5260, 5045] | Overflow | 6 |
| 5 | [5262, 5260], [5263, 5260], [5317, 5260], [5358, 5179], [5206, 5204] | Overflow | 4 |
| 5 | [7099, 7095],[7098, 7095],[7058, 7095], [6161, 7018],[7311, 7304] | Reserve Limit | NA |
| 5 | [7099, 7095],[7098, 7095],[7058, 7095], [7346, 7125],[7197, 7186] | Reserve Limit | NA |
| 5 | [5360, 5358], [8071, 8067], [5317, 5260], [5358, 5179],[5179, 5380] | Overflow | 1 |
| 6 | [5362, 5260], [5363, 5260], [5317, 5260], [5358, 5179],[5260, 5045],[5179, 5380] | Overflow | 6 |
| 6 | [5362, 5260], [5363, 5260], [5317, 5260], [5358, 5179],[5206, 5204],[5207, 5204] | Overflow | 4 |
| 6 | [7099, 7095], [7098, 7095], [7058, 7095], [6161, 7018],[7311, 7304],[7046, 7042] | Reserve Limit | NA |
| 6 | [7099, 7095], [7098, 7095], [7058, 7095], [7346, 7125],[7197, 7186],[7196, 7186] | Reserve Limit | NA |
| 7 | [5262, 5260], [5263, 5260], [5360, 5358], [5303, 5295], [5302, 5295], [5317, 5260], [5358, 5179] | Overflow<br>Undervoltage | 5<br>14 |
| 7 | [5262, 5260],[5263, 5260],[5317, 5260], [5358, 5179], [5260, 5045],[5179, 5380], [5380, 5384] | Overflow | 6 |
| 7 | [5262, 5260],[5263, 5260],[8071, 8067], [5068, 5063], [5067, 5063],[5069, 5063], [5317, 5260] | Overflow | 1 |
| 7 | [5262, 5260],[5263, 5260],[5317, 5260], [5358, 5179],[5206, 5204],[5207, 5204], [5208, 5204] | Overflow | 5 |
| 7 | [7099, 7095], [7098, 7095], [7058, 7095], [6161, 7018],[7311, 7304],[7046, 7042], [7045, 7042] | Reserve Limit | NA |
| 7 | [7099, 7095], [7098, 7095], [7058, 7095], [7311, 7304],[7046, 7042],[7045, 7042], [7304, 7095] | Reserve Limit | NA |
| 8 | [5262, 5260],[5263, 5260],[5360, 5358], [5303, 5295],[5302, 5295],[5317, 5260], [5358, 5179],[5260, 5045] | Overflow<br>Undervoltage | 5<br>14 |
| 8 | [5262, 5260],[5263, 5260],[5317, 5260], [5358, 5179],[5206, 5204],[5207, 5204], [5208, 5204],[5260, 5045] | Overflow | 5 |
| 8 | [5262, 5260],[5263, 5260],[5317, 5260], [5358, 5179],[5260, 5045],[5179, 5380], [5380, 5384],[5047, 5260] | Overflow | 6 |
| 8 | [5262, 5260],[5263, 5260],[8071, 8067], [5068, 5063], [5067, 5063],[5069, 5063], [5317, 5260],[5358, 5179] | Unsolved | NA |
| 8 | [6146, 6141],[6145, 6141],[6246, 6239], [6245, 6239],[6244, 6239],[6147, 6141], [6149, 6141],[6161, 7018] | Overflow | 1 |
| 8 | [7099, 7095], [7098, 7095], [7058, 7095], [6161, 7018],[7311, 7304],[7046, 7042], [7045, 7042],[6239, 7414] | Reserve Limit | NA |

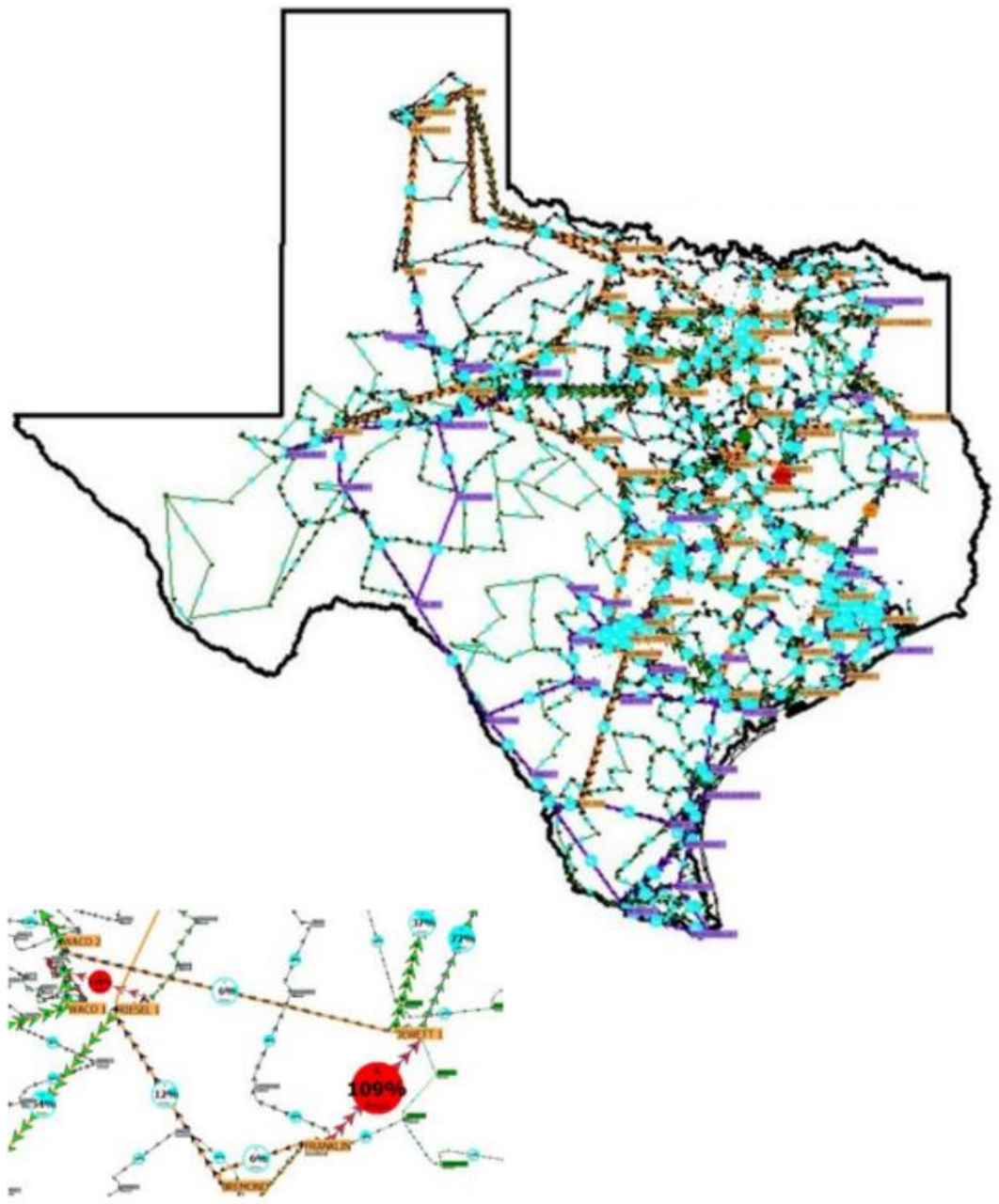

Figure 9. One-line diagram of the 2000-bus test cases and the corresponding overflow contingency caused by the removal of [5262, 5260],[5263, 5260],[5317, 5260],[5358, 5179] lines. The bottom portion of Figure 9 shows the zoom-in view of the area that violation has occurred.

## V. CONCLUSION

The proposed approach exploits the physical and topological characteristics of the electric power grids for identifying the critical lines in N-x contingency analysis. Actually, augmenting the LODF metric, that represent the physics of the electric grid, with the group betweenness centrality metric, that represent the topology of the electric grid, enables the proposed approach to effectively finds the critical lines even in large test cases. Protecting the resulting elements can enable power grid operators to prevent cascading failures and operate the system more reliably. The curse of dimensionality inherent to the generalized contingency analysis (i.e. N-x contingency analysis with $x \geq 2$) makes this problem computationally infeasible, for even midsize electric grid, to be solved by the traditional approaches. The proposed approach decouples the computation with the problem size which enables that to perform the N-x contingency analysis with $x \geq 2$ in a reasonable time irrespective of the grid size. Results show that the proposed approach acts as a straightforward and computationally tractable search engine that can quickly identify critical lines even in large test cases. Unfortunately, the computational infeasibility of traditional methods of contingency analysis makes the validation of the proposed approach challenging. In this connection, using multiple synthetic testcases and different contingency analysis, we demonstrated that our approach computes the realistic solutions and holds a promise for larger test cases. Our ongoing work aims to improve our approach by making it as independent as possible from "search level". In other words, we are working on an approach that can search in different part of electric power grids at the same time. Other ongoing work is developing further improvements to the proposed approach by updating the LODF metric after each line outage.

## REFERENCES

[1] V. Donde, V. Lopez, B. Lesieutre, A. Pinar, C. Yang, and J. Meza, "Severe multiple contingency screening in electric power systems," IEEE Transactions on Power Systems, vol. 23, no. 2, pp. 406–417, May 2008.

[2] X. Chen, K. Sun, Y. Cao, and S. Wang, "Identification of vulnerable lines in power grid based on complex network theory," in2007 I EE Power Engineering Society General Meeting, June 2007, pp. 1–6.

[3] B. C. Lesieutre, S. Roy, V. Donde, and A. Pinar, "Power system extreme event screening using graph partitioning," in2006 38th North American Power Symposium, 2006, pp. 503–510.

[4] E. P. R. Coelho, M. H. M. Paiva, M. E. V. Segatto, and G. Caporossi, "A new approach for contingency analysis based on centrality measures," IEEE Systems Journal, vol. 13, no. 2, pp. 1915–1923, June 2019.

[5] J. A. Kersulis, I. A. Hiskens, C. Coffrin, and D. K. Molzahn, "Topological graph metrics for detecting grid anomalies and improving algorithms," in2018 Power Systems Computation Conference (PSCC), June2018, pp. 1–7.

[6] S. Wasserman and K. Faust, "Social network analysis: methods and applications," Cambridge University Press, 1994.

[7] J. Yan, H. He, and Y. Sun, "Integrated security analysis on cascading failure in complex networks," IEEE Transactions on Information Forensics and Security, vol. 9, no. 3, pp. 451–463, March 2014.

[8] Z. Wang, A. Scaglione, and R. J. Thomas, "Electrical centrality measures for electric power grid vulnerability analysis," in49th IEEE Conference on Decision and Control (CDC), Dec 2010, pp. 5792–5797.

[9] P. Hines and S. Blumsack, "A centrality measure for electrical networks," in Proceedings of the 41st Annual Hawaii International Conference on System Sciences (HICSS 2008), Jan 2008, pp. 185–185.

[10] P. Hines, S. Blumsack, E. C. Sanchez, and C. Barrows, "The topological and electrical structure of power grids," in2010 43rd Hawaii International Conference on System Sciences, Jan 2010, pp. 1–10.

[11] A. Torres and G. Anders, "Spectral graph theory and network dependability," in2009 Fourth International Conference on Dependability of Computer Systems, June 2009, pp. 356–363.

[12] Y. C. B. K. S. J. D. C. D. B. I. Gorton, Z. Huang and J. Feo, "A high-performance hybrid computing approach to massive contingency analysis in the power grid," December 2009, p. 277–283.

[13] M. Halappanavar, Y. Chen, R. Adolf, D. Haglin, Z. Huang, and M. Rice, "Towards efficient n-x contingency selection using group betweenness centrality," pp. 273–282, Nov 2012.

[14] H. Bai and S. Miao, "An electrical betweenness approach for vulnerability assessment of power grids considering the capacity of generators and load," IET Generation, Transmission Distribution, vol. 9, no. 12, pp. 1324–1331, 2015.

[15] K. Wang, B. Zhang, Z. Zhang, X. Yin, and B. Wang, "Hybrid flow betweenness approach for identification of vulnerable line in power system," Physica A, vol. 390, no. 23, pp. 4692–4701, 2011.

[16] E. Bompard, D. Wu, and F. Xue, "Structural vulnerability of power systems: A topological approach," Electric Power Systems Research, vol. 81, no. 7, pp. 1334 – 1340, 2011. [Online]. Available:http://www.sciencedirect.com/science/article/pii/S0378779611000332

[17] G. chen, Z. Yang Dong, D. J. Hill, and G. Hua Zhang, "An improved model for structural vulnerability analysis of power networks," Physica A, vol. 388, no. 23, pp. 4259–4266, 2009.

[18] C. M. Davis and T. J. Overbye, "Multiple element contingency screening," IEEE Transactions on Power Systems, vol. 26, no. 3, pp. 1294–1301, Aug 2011.

[19] M. Everett and S. Borgatti, "The centrality of groups and classes," Journal of Mathematical Sociology, vol. 23, no. 9, pp. 181–201, 1999.

[20] S. Borgatti, "Identifying sets of key players in a social network," Computational & Mathematical Organization Theory, vol. 12, no. 1,pp. 21–34, 2006.

[21] PowerWorld Corporation, available at http://www.powerworld.com,2018.

[22] T. Guler and G. Gross, "Detection of island formation and identification of causal factors under multiple line outages," IEEE Transactions on Power Systems, vol. 22, no. 2, pp. 505–513, May 2007.
[23] D. A. Tejada-Arango, P. Sánchez-Martın, and A. Ramos, "Security constrained unit commitment using line outage distribution factors,"IEEE Transactions on Power Systems, vol. 33, no. 1, pp. 329–337, Jan2018.
[24] M. O. W. Grond, J. I. P. Pouw, J. Morren, and H. J. G. Slootweg,"Applicability of line outage distribution factors to evaluate distribution network expansion options," in2014 49th International Universities Power Engineering Conference (UPEC), Sep. 2014, pp. 1–6.
[25] R. Baldick, "Variation of distribution factors with loading," IEEE Transactions on Power Systems, vol. 18, no. 4, pp. 1316–1323, Nov2003.
[26] A. Al-Digs, S. V. Dhople, and Y. C. Chen, "Dynamic distribution factors," IEEE Transactions on Power Systems, vol. 34, no. 6, pp. 4974–4983, Nov 2019.
[27] J. Guo, Y. Fu, Z. Li, and M. Shahidehpour, "Direct calculation of line outage distribution factors," IEEE Transactions on Power Systems, vol. 24, no. 3, pp. 1633–1634, Aug 2009.
[28] P. Hines, I. Dobson, and P. Rezaei, "Cascading power outages propagate locally in an influence graph that is not the actual grid topology," in2017IEEE Power Energy Society General Meeting, July 2017, pp. 1–1.
[29] benchmark library for electric power grids in Texas A&M University, available at https://electricgrids.engr.tamu.edu.
[30] A. B. Birchfield, T. Xu, K. M. Gegner, K. S. Shetye, and T. J. Over-bye, "Grid structural characteristics as validation criteria for synthetic networks," IEEE Transactions on Power Systems, vol. 32, no. 4, pp.3258–3265, July 2017.
[31] B. L. Thayer, Z. Mao, Y. Liu, K. Davis, and T. Overbye, "Easy simauto (esa): A python package that simplifies interacting with powerworld simulator," Journal of Open Source Software, vol. 5, no. 50, p. 2289,2020. [Online]. Available: https://doi.org/10.21105/joss.02289.